\documentclass[12pt]{article}

\usepackage{amsmath,amsfonts,amssymb,amsthm,amsxtra} 
\usepackage{amssymb}
\usepackage{hyperref}

\makeatletter
\newfont{\footsc}{cmcsc10 at 8truept}
\newfont{\footbf}{cmbx10 at 8truept}  
\newfont{\footrm}{cmr10 at 10truept} 
\makeatother
\def\stf#1#2{\left[#1\atop#2\right]} 
\def\sts#1#2{\left\{#1\atop#2\right\}}

\newtheorem{theorem}{Theorem}
\newtheorem{Prop}{Proposition}

\newtheorem{Lem}{Lemma}

\begin{document}

\title{Cauchy-Carlitz numbers}
\author{
Hajime Kaneko\\ 
\small Institute of Mathematics\\
\small University of Tsukuba\\
\small 1-1-1 Tennodai, Tsukuba, Ibaraki 350-0006 JAPAN\\
\small Center for Integrated Research in Fundamental Science and Engineering (CiRfSE)\\
\small University of Tsukuba, \\
\small Tsukuba, Ibaraki 305-8571, JAPAN\\
%\small Center for Integrated Research 
%in Fundamental Science and Technology (CiRfSE) \\
%\small University of Tsukuba \\
%\small Tsukuba, Ibaraki 305-8571 JAPAN\\
\small \texttt{kanekoha@math.tsukuba.ac.jp}\\\\
Takao Komatsu\\ 
\small School of Mathematics and Statistics\\
\small Wuhan University\\
\small Wuhan 430072 China\\
\small \texttt{komatsu@whu.edu.cn}
}

\date{
%\small Submitted: November 1, 2014;  Accepted: December 2, 2014.\\
%\small MR Subject Classifications: Primary 11B73; Secondary 11B68, 11R58
}

\maketitle

\begin{abstract}
In 1935 Carlitz introduced Bernoulli-Carlitz numbers as analogues of Bernoulli numbers for the rational function field $\mathbb F_r(T)$.  In this paper, we introduce Cauchy-Carlitz numbers as analogues of Cauchy numbers. By using Stirling-Carlitz numbers, we give their arithmetical and combinatorial properties and relations with Bernoulli-Carlitz numbers for $\mathbb F_r(T)$. Several new identities are also obtained by using Hasse-Teichim\"uller derivatives. 

\end{abstract}

\section{Introduction}

In 1935, L. Carlitz (\cite{Car35}) introduced analogues of Bernoulli numbers for the rational function field $K=\mathbb F_r(T)$, which are called Bernoulli-Carlitz numbers now. He proved an analogue of the von Staudt-Clausen theorem (\cite{Car37,Car40}). Some identities for Bernoulli-Carlitz numbers were found in \cite{Gekeler}. In \cite{JKS}, explicit formulae of Bernoulli-Carlitz numbers were given by using the basic properties of the Hasse-Teichm\"uller derivatives. In \cite{Taelman}, it was shown a necessary and sufficient condition for
a nonzero prime ideal of the rational function field divides the $n$-th Bernoulli-Carlitz number. A recent exposition of Bernoulli-Carlitz numbers can be seen in \cite{Rodriguez}. We refer to \cite{Goss} for an exposition and the modern notation.  
The Carlitz exponential $e_C(x)$ is defined by 
\begin{equation} 
e_C(x)=\sum_{i=0}^\infty\frac{x^{r^i}}{D_i}\,, 
\label{carlitzexp}
\end{equation}  
where $D_i=[i][i-1]^r\cdots [1]^{r^{i-1}}$ ($i\ge 1$) with $D_0=1$, and $[i]=T^{r^i}-T$.  
The Carlitz logarithm $\log_C(x)$ is defined by 
\begin{equation} 
\log_C(x)=\sum_{i=0}^\infty(-1)^i\frac{x^{r^i}}{L_i}\,, 
\label{carlitzlog}
\end{equation}  
where $L_i=[i][i-1]\cdots [1]$ ($i\ge 1$) with $L_0=1$.  
Notice that 
\begin{equation}  
e_C\bigl(\log_C(x)\bigr)=\log_C\bigl(e_C(x)\bigr)=x\,. 
\label{rel:explog}
\end{equation}  
$e_C(x)$ and $\log_C(x)$ have the functional equations 
\begin{equation} 
e_C(Tx)=T e_C(x)+e_C(x)^r 
\label{funceq:exp} 
\end{equation}  
and 
\begin{equation} 
T\log_C(x)=\log_C(Tx)+\log_C(x^r)\,, 
\label{funceq:exp} 
\end{equation}  
respectively.  

The Carlitz factorial $\Pi(i)$ is defined by 
\begin{equation}  
\Pi(i)=\prod_{j=0}^m D_j^{c_j}
\label{carlitzfac}
\end{equation} 
for a non-negative integer $i$ with $r$-ary expansion: 
\begin{equation}  
i=\sum_{j=0}^m c_j r^j\quad (0\le c_j<r)\,. 
\label{r-expansion}
\end{equation}  
Therefore, 
\begin{align}
\Pi(i)&=\prod_{j=0}^m\left(\prod_{k=1}^j(T^{r^j}-T^{r^{j-k}})\right)^{c_j}\notag\\
&=\prod_{k=1}^m(T^{r^k}-T)^{c_k+c_{k+1}r+\cdots}
=\prod_{k=1}^m(T^{r^k}-T)^{[i/r^k]}\quad\hbox{($[\cdot]$ denotes the floor function.)}\notag\\
&=\prod_{k=1}^m(T^{r^k}-T)^{c_k+c_{k+1}r+\cdots+c_m r^{m-k}}   
\label{carlitzfac2} 
\end{align} 
and 
\begin{equation} 
\Pi(r^d-1)=(D_0\cdots D_{d-1})^{r-1}=\frac{D_d}{L_d}\quad (d\ge 0)\,. 
\label{rel:dl}
\end{equation}

\section{Cauchy-Carlitz numbers}   

The Bernoulli-Carlitz numbers $BC_n$ are defined by 
\begin{equation}  
\frac{x}{e_C(x)}=\sum_{n=0}^\infty\frac{BC_n}{\Pi(n)}x^n
\label{bercarlitz} 
\end{equation}  
as analogues of the classical Bernoulli numbers $B_n$, defined by 
$$
\frac{x}{e^x-1}=\sum_{n=0}^\infty\frac{B_n}{n!}x^n\,. 
$$ 

As analogues of the classical Cauchy numbers $c_n$, defined by 
\begin{equation} 
\frac{x}{\log(1+x)}=\sum_{n=0}^\infty\frac{c_n}{n!}x^n\,, 
\label{classicalcauchyfirst} 
\end{equation}  
we define the Cauchy-Carlitz numbers $CC_n$ by 
\begin{equation}  
\frac{x}{\log_C(x)}=\sum_{n=0}^\infty\frac{CC_n}{\Pi(n)}x^n\,. 
\label{caucarlitz} 
\end{equation}  

In addition, as analogues of the Stirling numbers of the first kind $\stf{n}{k}$, defined by 
\begin{equation} 
\frac{\bigl(-\log(1-t)\bigr)^k}{k!}=\sum_{n=0}^\infty\stf{n}{k}\frac{t^n}{n!}\,,
\label{stf}  
\end{equation}
we define the Stirling-Carlitz numbers of the first kind  $\stf{n}{k}_C$ by 
\begin{equation}  
\frac{\bigl(\log_C(z)\bigr)^k}{\Pi(k)}=\sum_{n=0}^\infty\stf{n}{k}_C\frac{z^n}{\Pi(n)}\,. 
\label{stfcarlitz}  
\end{equation}  
As analogues of the Stirling numbers of the second kind $\sts{n}{k}$, defined by
$$
\frac{(e^t-1)^k}{k!}=\sum_{n=0}^\infty\sts{n}{k}\frac{t^n}{n!}\,,
$$ 
we define the Stirling-Carlitz numbers of the second kind  $\sts{n}{k}_C$ by 
\begin{equation}  
\frac{\bigl(e_C(z)\bigr)^k}{\Pi(k)}=\sum_{n=0}^\infty\sts{n}{k}_C\frac{z^n}{\Pi(n)}\,. 
\label{stscarlitz}  
\end{equation}  
By the definition (\ref{stfcarlitz}), we have 
\begin{equation}
\stf{n}{0}_C=0\quad(n\ge 1),\quad \stf{n}{m}_C=0\quad(n<m)\quad  
\hbox{and}\quad \stf{n}{n}_C=1\quad(n\ge 0)\,.
\label{eqn:aabb} 
\end{equation}
Similarly, we see 
\begin{equation}
\sts{n}{0}_C=0\quad(n\ge 1),\quad \sts{n}{m}_C=0\quad(n<m)\quad  
\hbox{and}\quad \sts{n}{n}_C=1\quad(n\ge 0)\,.
\label{eqn:aacc} 
\end{equation}
It is known that poly-Cauchy numbers $\mathfrak c_n^{(k)}$ 
are expressed in terms of the Stirling numbers of the first kind:  
\begin{equation}
\mathfrak c_n^{(k)}=\sum_{m=0}^n\stf{n}{m}\frac{(-1)^{n-m}}{(m+1)^k}
\label{eqn:aadd}
\end{equation}
(\cite[Theorem 1]{Ko1})\footnote{In this paper we use the notation $\mathfrak c_n^{(k)}$ in order to distinguish from the Cauchy numbers of high order in later section.}.  If $k=1$, this is an explicit expression of the classical Cauchy numbers $c_n$ (\cite[Ch. VII]{Comtet},\cite[p.1908]{MSV}).

As analogues, we express Cauchy-Carlitz numbers as certain finite sums of 
Stirling-Carlitz numbers of the first kind. 
Thus, Cauchy-Carlitz numbers are calculated by this expression because 
Stirling-Carlitz numbers of the first kind are computed by (\ref{stfcarlitz}). 

\begin{theorem} 
\begin{equation}  
CC_n=\sum_{j=0}^\infty\frac{1}{L_j}\stf{n}{r^j-1}_C\,.  
\label{carlitzcauchystf}
\end{equation} 
\label{th10}
\end{theorem}  
\begin{proof}  
Note that the right-hand side of (\ref{carlitzcauchystf}) is a finite sum by the second relation of 
(\ref{eqn:aabb}). Observe that 
\begin{align*} 
\frac{z}{\log_C(z)}&=\frac{e_C\bigl(\log_C(z)\bigr)}{\log_C(z)}=\left.\frac{e_C(t)}{t}\right|_{t=\log_C(z)}\\
&=\sum_{j=0}^\infty\frac{\bigl(\log_C(z)\bigr)^{r^j-1}}{D_j}\\
&=\sum_{j=0}^\infty\frac{1}{D_j}\Pi(r^j-1)\sum_{n=0}^\infty\stf{n}{r^j-1}_C\frac{z^n}{\Pi(n)}\\
&=\sum_{j=0}^\infty\frac{1}{L_j}\sum_{n=0}^\infty\stf{n}{r^j-1}_C\frac{z^n}{\Pi(n)}\\ 
&=\sum_{n=0}^\infty\left(\sum_{j=0}^\infty\frac{1}{L_j}\stf{n}{r^j-1}_C\right)\frac{z^n}{\Pi(n)}\,.  
\end{align*} 
By the definition (\ref{caucarlitz}), we get (\ref{carlitzcauchystf}).  
\end{proof}

Similarly, as an analogous expression of the classical Bernoulli numbers
$$
B_n=\sum_{m=0}^n\sts{n}{m}\frac{(-1)^{n-m}m!}{m+1}\,,
$$  
%we have the following.
we get that Bernoulli-Carlitz numbers are equal to certain finite sums of 
Stirling-Carlitz numbers of the second kind. In particular, Bernoulli-Carlitz numbers 
are calculated by (\ref{stscarlitz}).   

\begin{theorem} 
\begin{equation}  
BC_n=\sum_{j=0}^\infty\frac{(-1)^j D_j}{L_j^2}\sts{n}{r^j-1}_C\,.  
\label{carlitzbernoullistf} 
\end{equation} 
\label{th20}
\end{theorem}

Let $f(z)$ be a formal power series of the form $f(z)=\sum_{i=0}^{\infty} f_i z^{r^i}\in \mathbb F_r(T)[[z]]$. 
If $f_0\ne 0$, then let $g(z)=\sum_{i=0}^{\infty}g_i z^{r^i}$ be the inverse function of $f(z)$. 
Set 
\begin{align}
h(z)=\sum_{n=0}^{\infty} h_n z^n:=\frac{z f'(z)}{f(z)}\,.
\label{general}
\end{align}
Then (\ref{bercarlitz}) and (\ref{caucarlitz}) are special examples of (\ref{general}). 
Carlitz \cite{Car35} studied the coefficients $h_n$ in the case where $n$ satisfies certain assumptions. 
For instance, he showed that if $f_0=1$, then $h_{r^k-1}=g_k$. Similarly, it is seen that if $f_0\ne 0$, then 
$h_{r^k-1}=f_0^{r^k} g_k$. 

J. A. Lara Rodr\'iguez and D. S. Thakur \cite{RT} gave other relations as follows: 
Let $l$ be a integer with $1\leq l\leq r$ and $k,k_1,\ldots,k_l$ integers with $0\leq k_j \leq k$ 
for any $1\leq j\leq l$. Then 
$$
\prod_{j=1}^l h_{r^k-r^{k_j}}=h_{\sum (r^k-r^{k_j})} \,.
$$

\section{Examples}  

In this section we give examples of Stirling-Carlitz numbers of the first and second kind. 
Moreover, by using Theorem \ref{th10} and Theorem \ref{th20}, we calculate examples of 
Cauchy-Carlitz numbers and Bernoulli-Carlitz numbers. 
In the rest of this section, we assume that $r=3$. \par
Observe that 
\begin{align*}
\sum_{n=0}^{\infty}\stf{n}{2}_C \frac{z^n}{\Pi(n)} & = \frac{\bigl(\log_C(z)\bigr)^2}{\Pi(2)}\\
& = \left(z-\frac{1}{L_1} z^3+\frac{1}{L_2} z^9-+\cdots\right)^2\\
& = z^2-\frac{2}{[1]}z^4+\frac{1}{[1]^2}z^6+0\cdot z^8+\cdots\,. 
\end{align*}
Hence, we get 
\begin{align*}
\stf{4}{2}_C=-\frac{2}{[1]}\Pi(4)=-2=1, \quad 
\stf{6}{2}_C=\frac{1}{[1]^2}\Pi(6)=1\quad \hbox{and}\quad \stf{8}{2}_C=0\,. 
\end{align*}
By combining Theorem \ref{th10}, 
(\ref{eqn:aabb}), and the equality above, we obtain 
$$
CC_2=\frac{1}{L_0}\stf{2}{0}_C+\frac{1}{L_1}\stf{2}{2}_C=\frac{1}{T^3-T}=\frac{1}{T^3+2 T}\,, 
$$ 
$$
CC_4=\frac{1}{L_1}\stf{4}{2}_C=\frac{1}{T^3+2 T}\,, \ 
CC_6=\frac{1}{L_1}\stf{6}{2}_C=\frac{1}{T^3+2 T} \,,
$$ 
and 
\begin{align*} 
CC_8&=\frac{1}{L_1}\stf{8}{2}_C+\frac{1}{L_2}\stf{8}{8}_C=\frac{1}{L_2}\\
&=\frac{1}{(T^3-T)(T^9-T)}=\frac{1}{(T^3+2 T)(T^9+2 T)}\,. 
\end{align*}  
In the same way, by using 
\begin{align*}
\sum_{n=0}^{\infty}\sts{n}{2}_C \frac{z^n}{\Pi(n)} & = \frac{\bigl(e_C(z)\bigr)^2}{\Pi(2)}
= z^2+\frac{2}{[1]}z^4+\frac{1}{[1]^2}z^6+0\cdot z^8+\cdots, 
\end{align*}
we get 
\begin{align*}
\sts{4}{2}_C=2, \quad 
\sts{6}{2}_C=1\quad \hbox{and}\quad \sts{8}{2}_C=0\,. 
\end{align*}
By combining Theorem \ref{th20}, (\ref{eqn:aacc}), and the equality above, we obtain the following: 
\begin{align*} 
BC_2&=-\frac{D_1}{L_1^2}\sts{2}{2}_C=\frac{-1}{T^3-T}=\frac{2}{T^3+2 T}\,,\\
BC_4&=-\frac{D_1}{L_1^2}\sts{4}{2}_C=\frac{1}{T^3+2 T}\,, \ 
BC_6=-\frac{D_1}{L_1^2}\sts{6}{2}_C=\frac{2}{T^3+2 T}\,,\\
BC_8&=-\frac{D_1}{L_1^2}\sts{8}{2}_C+\frac{D_2}{L_2^2}\sts{8}{8}_C=\frac{T^3-T}{T^9-T}=\frac{1}{T^6+T^4+T^2+1}\,.\\
\end{align*}

\section{Hasse-Teichm\"uller derivatives}  

%Let $\mathbb F$ be a field, $\mathbb F[z]$ be a polynomial ring in one variable $z$ over $\mathbb F$. 
%For each integer $n\ge 0$, define the Hasse-Teichm\"uller derivative $H^{(n)}$ of order $n$ by 
%$$
%H^{(n)}(z^m)=\binom{m}{n}z^{m-n} 
%$$  
%for any integer $m\ge 0$.  
Let $\mathbb{F}$ be a field of any characteristic, $\mathbb{F}[[z]]$ the ring of formal power series 
in one variable $z$, and $\mathbb{F}((z))$ the field of Laurent series in $z$. Let $n$ be a nonnegative 
integer. We define the Hasse-Teichm\"uller derivative $H^{(n)}$ of order $n$ by 
$$
H^{(n)}\left(\sum_{m=R}^{\infty} c_m z^m\right)
=\sum_{m=R}^{\infty} c_m \binom{m}{n}z^{m-n}
$$
for $\sum_{m=R}^{\infty} c_m z^m\in \mathbb{F}((z))$, 
where $R$ is an integer and $c_m\in\mathbb{F}$ for any $m\geq R$. 

The Hasse-Teichm\"uller derivatives satisfy the product rule \cite{Teich}, the quotient rule \cite{GN} and the chain rule \cite{Hasse}. 
One of the product rules can be described as follows.  
\begin{Lem}  
For $f_i\in\mathbb F[[z]]$ ($i=1,\dots,k$) with $k\ge 2$ and for $n\ge 1$, we have 
$$
H^{(n)}(f_1\cdots f_k)=\sum_{i_1,\dots,i_k\ge 0\atop i_1+\cdots+i_k=n}H^{(i_1)}(f_1)\cdots H^{(i_k)}(f_k)\,. 
$$ 
\label{productrule2}
\end{Lem} 

The quotient rules can be described as follows.  
\begin{Lem}  
For $f\in\mathbb F[[z]]\backslash \{0\}$ and $n\ge 1$,  
we have 
\begin{align} 
H^{(n)}\left(\frac{1}{f}\right)&=\sum_{k=1}^n\frac{(-1)^k}{f^{k+1}}\sum_{i_1,\dots,i_k\ge 1\atop i_1+\cdots+i_k=n}H^{(i_1)}(f)\cdots H^{(i_k)}(f)
\label{quotientrule1}\\ 
&=\sum_{k=1}^n\binom{n+1}{k+1}\frac{(-1)^k}{f^{k+1}}\sum_{i_1,\dots,i_k\ge 0\atop i_1+\cdots+i_k=n}H^{(i_1)}(f)\cdots H^{(i_k)}(f)\,.
\label{quotientrule2} 
\end{align}   
\label{quotientrules}
\end{Lem}  

By using the Hasse-Teichm\"uller derivative of order $n$, some explicit expressions of Bernoulli-Carlitz numbers are obtained in \cite{JKS}. In this section %we obtain another explicit expression of Cauchy-Carlitz numbers. 
we study explicit expressions of Cauchy-Carlitz numbers and Cauchy numbers, by using 
the Hasse-Teichm\"uller derivatives. First, we express Cauchy-Carlitz numbers 
in terms of \(L_i\) without using Stirling-Carlitz numbers of the first kind, 
which gives a new method to calculate Cauchy-Carlitz numbers. 
%We apply the method to \(CC_8\) in Section 3. 

\begin{theorem}  
For $n\ge 1$, 
$$
CC_n=\Pi(n)\sum_{k=1}^n(-1)^k\sum_{i_1,\ldots,i_k\geq 1%r^{i_1},\dots,r^{i_k}\ge 2
\atop r^{i_1}+\cdots+r^{i_k}=n+k}\frac{(-1)^{i_1+\cdots+i_k}}{L_{i_1}\cdots L_{i_k}}\,.
$$ 
\label{cc_ht}
\end{theorem}  
\begin{proof}  
Put 
$$
g:=\frac{\log_C(z)}{z}=\sum_{j=0}^\infty(-1)^j\frac{z^{r^j-1}}{L_j}\,. 
$$ 
Note that 
\begin{align} 
\left.H^{(e)}(g)\right|_{z=0}&=\left.\sum_{j=0}^\infty\frac{(-1)^j}{L_j}\binom{r^j-1}{e}z^{r^j-1-e}\right|_{z=0}
\nonumber \\
&=\left\{
\begin{array}{cc}
%\frac{(-1)^i}{L_i} 
\frac{(-1)^i}{L_i} & \hbox{if }e=r^i-1\,, \\
0 & \hbox{otherwise}\,.
\end{array}
\right.
\label{coefficientg}
\end{align}  
%where $e=r^i-1$.   
Thus, by %applying Lemma \ref{productrule2} to (\ref{caucarlitz}) and 
using Lemma \ref{quotientrules} (\ref{quotientrule1}) and (\ref{coefficientg}), we have 
\begin{align*} 
\frac{CC_n}{\Pi(n)}&=\left.H^{(n)}\left(\frac{z}{\log_C(z)}\right)\right|_{z=0}=\left.H^{(n)}\left(\frac{1}{g}\right)\right|_{z=0}\\
&=\sum_{k=1}^n(-1)^k\sum_{e_1,\dots,e_k\ge 1\atop e_1+\cdots+e_k=n}\left.H^{(e_1)}(g)\right|_{z=0}\cdots\left.H^{(e_k)}(g)\right|_{z=0}\\
&=\sum_{k=1}^n(-1)^k\sum_{i_1,\ldots,i_k\geq 1%r^{i_1},\dots,r^{i_k}\ge 2
\atop r^{i_1}+\cdots+r^{i_k}=n+k}\frac{(-1)^{i_1+\cdots+i_k}}{L_1\cdots L_{i_k}}\,.  
\end{align*} 
\end{proof}  

\noindent 
{\bf Example.}   
Let $r=3$ and $n=8$.  
%Since $i_1=2$ for $k=1$ and $i_1=i_2=i_3=i_4=1$ for $k=4$, we have 
For $1\leq k\leq n$, put 
$$
S_k=\{
(i_1,\ldots,i_k)\mid i_1,\ldots,i_k\geq 1, \ 3^{i_1}+\cdots+3^{i_k}=8+k
\}\,.
$$
Then $S_k$ is empty except the cases of $k=1, S_1=\{(2)\}$ and 
$k=4, S_4=\{(1,1,1,1)\}$. By Theorem \ref{cc_ht}, we get 
\begin{align*} 
CC_8&=\Pi(8)\left((-1)^1\frac{(-1)^2}{L_2}+(-1)^4\frac{(-1)^4}{L_1^4}\right)\\
&=(T^3-T)^2\left(-\frac{1}{(T^3-T)(T^9-T)}+\frac{1}{(T^3-T)^4}\right)\\
&=\frac{1}{(T^3-T)(T^9-T)}\,. 
\end{align*}

%The method to obtain an explicit formula for Cauchy-Carlitz numbers can be applicable to that of the classical Cauchy numbers.  
%It is known that Cauchy numbers (of the first kind) can be expressed in terms of the (unsigned) Stirling numbers of the first kind:  
%$$
%c_n=\sum_{m=0}^n\stf{n}{m}\frac{(-1)^{n-m}}{m+1}\quad(n\ge 0)
%$$ 
%(\cite[Ch. VII]{Comtet},\cite[Theorem 1]{Ko1},\cite[p.1908]{MSV}). 

%By using the Hasse-Teichm\"uller derivative of order $n$, we obtain another explicit expression of classical Cauchy numbers. 

Recall that if $k=1$, then (\ref{eqn:aadd}) gives an explicit formula for the classical Cauchy numbers 
in terms of the Stirling numbers of the first kind. 
By applying the method for the proof of Theorem \ref{cc_ht}, we get an explicit formula for the 
classical Cauchy numbers without using the Stirling numbers of the first kind. 

\begin{theorem}  
For $n\ge 1$, 
$$
c_n=(-1)^n n!\sum_{k=1}^n(-1)^k\sum_{i_1,\dots,i_k\ge 2\atop i_1+\cdots+i_k=n+k}\frac{1}{i_1\cdots i_k}\,.
$$ 
\label{th_classicalcauchy1}
\end{theorem}  
\begin{proof}  
Put 
$$
h:=\frac{\log(1+z)}{z}=\sum_{j=0}^\infty(-1)^j\frac{z^j}{j+1}\,. 
$$ 
Note that 
\begin{align*} 
\left. H^{(i)}(h)\right|_{z=0}=\left.\sum_{j=0}^\infty\frac{(-1)^j}{j+1}\binom{j}{i}z^{j-i}\right|_{z=0}
=\frac{(-1)^i}{i+1}\,.
\end{align*}    
Hence, by using Lemma \ref{quotientrules} (\ref{quotientrule1}), we have 
\begin{align*} 
\frac{c_n}{n!}&=\left.H^{(n)}\left(\frac{z}{\log(1+z)}\right)\right|_{z=0}=\left.H^{(n)}\left(\frac{1}{h}\right)\right|_{z=0}\\
&=\sum_{k=1}^n(-1)^k\sum_{i_1,\dots,i_k\ge 1\atop i_1+\cdots+i_k=n}\left.H^{(i_1)}(h)\right|_{z=0}\cdots\left.H^{(i_k)}(h)\right|_{z=0}\\
&=\sum_{k=1}^n(-1)^k\sum_{i_1,\dots,i_k\ge 1\atop i_1+\cdots+i_k=n}\frac{(-1)^{i_1+\cdots+i_k}}{(i_1+1)\cdots(i_k+1)}\\
&=(-1)^n\sum_{k=1}^n(-1)^k\sum_{i_1,\dots,i_k\ge 2\atop i_1+\cdots+i_k=n+k}\frac{1}{i_1\cdots i_k}\,. 
\end{align*} 
\end{proof}  

We can express the Cauchy numbers also in terms of the binomial coefficients. 
In fact, by using Lemma \ref{quotientrules} (\ref{quotientrule2}) 
instead of Lemma \ref{quotientrules} (\ref{quotientrule1}) 
in the proof of Theorem \ref{th_classicalcauchy1}, we obtain the following:  

\begin{Prop}  
For $n\ge 1$, 
$$
c_n=(-1)^n n!\sum_{k=1}^n(-1)^k\binom{n+1}{k+1}\sum_{i_1,\dots,i_k\ge 1\atop i_1+\cdots+i_k=n+k}\frac{1}{i_1\cdots i_k}\,.
$$ 
\label{th_classicalcauchy2}
\end{Prop}  

From (\ref{stf}) we have  
\begin{equation}  
\left(\frac{-\log(1-z)}{z}\right)^k=\sum_{n=k}^\infty k!\stf{n}{k}\frac{z^{n-k}}{n!}
=\sum_{n=0}^{\infty} \frac{k !}{(n+k) !} \stf{n+k}{k} z^n\,. 
\label{eq19} 
\end{equation}  
Applying Lemma \ref{productrule2} with 
$$
f_1(z)=\cdots=f_k(z)=\frac{-\log(1-z)}{z}\,, 
$$
we get 
\begin{equation}  
\sum_{i_1,\dots,i_k\ge 0\atop i_1+\cdots+i_k=n}\frac{1}{(i_1+1)\cdots(i_k+1)}=\frac{k!}{(n+k)!}\stf{n+k}{k}\,. 
\label{eq20} 
\end{equation}   
Together with Proposition \ref{th_classicalcauchy2}, 
we deduce a simple expression for Cauchy numbers %of the first kind 
in terms of the binomial coefficients and Stirling numbers of the first kind.  

\begin{Prop} 
For $n\ge 1$ 
$$
c_n=(-1)^n\sum_{k=1}^n(-1)^k\dfrac{\binom{n+1}{k+1}}{\binom{n+k}{n}}\stf{n+k}{k}\,.
$$ 
\label{th_classicalcauchy3} 
\end{Prop}  

By using Proposition \ref{th_classicalcauchy3}, immediately we get some initial values of Cauchy numbers: 
$$
c_1=\frac{1}{2},~ c_2=-\frac{1}{6},~ c_3=\frac{1}{4},~ c_4=-\frac{19}{30},~ c_5=\frac{9}{4},~ c_6=-\frac{863}{84},~ c_7=\frac{1375}{24}. 
$$  
\bigskip 

Define the Cauchy numbers $c_n^{(m)}$ of order $m$ by 
\begin{equation}  
\left(\frac{z}{\log(1+z)}\right)^m=\sum_{n=0}^\infty c_n^{(m)}\frac{z^n}{n!}\,. 
\label{cauchyhighorder}
\end{equation}   
Notice that the concept of Cauchy numbers of higher order is different from that of poly-Cauchy numbers (\cite{Ko1}), though we use the similar notation here.  

%By applying Lemma \ref{quotientrules} (\ref{quotientrule1}) first and then Lemma \ref{productrule2} to (\ref{cauchyhighorder}), and using the identity (\ref{eq20}), we get an explicit formula for Cauchy numbers of higher order.  
We now introduce formulae for Cauchy numbers of higher order, 
by using Hasse-Teichm\"uller derivatives. The first two formulae (Proposition \ref{Proposition3} 
and Proposition \ref{Proposition4}) give expressions for Cauchy numbers in terms of 
multinomial coefficients, binomial coefficients, and the Stirling numbers of the first kind. 
The last formula (Proposition \ref{Proposition5}) gives a simple expression 
without using multinomial coefficients, which is useful to calculate Cauchy numbers of higher order. \par
Let again $h(z)=(\log(1+z))/z$. By applying Lemma \ref{quotientrules} (\ref{quotientrule1}) 
with \(f(z)=h(z)^m\), we see, by (\ref{cauchyhighorder}), 
$$
\frac{c_n^{(m)}}{n !}=\sum_{k=1}^n (-1)^k \sum_{i_1,\ldots,i_k\geq 1\atop i_1+\cdots+i_k=n}
\left. H^{(i_1)} \left( h^m\right)\right|_{z=0} \cdots \left. H^{(i_k)} \left( h^m\right)\right|_{z=0}.
$$
Using Lemma \ref{productrule2} and the identity (\ref{eq20}), 
we get an explicit formula for Cauchy numbers of higher order.  

\begin{Prop}  
For $n\ge 1$ 
$$
c_n^{(m)}=(-1)^n\sum_{k=1}^n(-1)^k\sum_{i_1,\dots,i_k\ge 1\atop i_1+\cdots+i_k=n}\dfrac{\binom{n}{i_1,\dots,i_k}}{\binom{i_1+m}{m}\cdots\binom{i_k+m}{m}}\stf{i_1+m}{m}\cdots\stf{i_k+m}{m}\,. 
$$ 
\label{Proposition3}
\end{Prop}

If we use Lemma \ref{quotientrules} (\ref{quotientrule2}) instead of Lemma \ref{quotientrules} (\ref{quotientrule1}), we obtain the following.  

\begin{Prop}  
For $n\ge 1$ 
$$ 
c_n^{(m)}=(-1)^n\sum_{k=1}^n(-1)^k\binom{n+1}{k+1}\sum_{i_1,\dots,i_k\ge 0\atop i_1+\cdots+i_k=n}\dfrac{\binom{n}{i_1,\dots,i_k}}{\binom{i_1+m}{m}\cdots\binom{i_k+m}{m}}\stf{i_1+m}{m}\cdots\stf{i_k+m}{m}\,. 
$$ 
\label{Proposition4}
\end{Prop}

%If we apply Lemma \ref{quotientrules} (\ref{quotientrule1}) to the expansion in (\ref{eq19}) and replacing $k$ by $m k$, we have a different explicit expression of $c_n^{(m)}$.  
Applying Lemma \ref{productrule2} with 
$$
f_1(z)=\cdots=f_k(z)=\left(-\frac{\log (1-z)}{z}\right)^m,
$$
we get, by (\ref{eq19}), 
$$
\frac{(mk)!}{(n+mk)!}\stf{n+mk}{mk} 
=
\sum_{i_1,\ldots,i_k\geq 0\atop i_1+\cdots+i_k=n} \frac{m!}{(i_1+m)!}\cdots\frac{m!}{(i_k+m)!}
\stf{i_1+m}{m}\cdots\stf{i_k+m}{m}\,. 
$$
Multiplying the both sides of the equality above by \(n!\), we deduce a different explicit expression of $c_n^{(m)}$ 
by Proposition \ref{Proposition4}. 

\begin{Prop}  
For $n\ge 1$ 
$$
c_n^{(m)}=(-1)^n\sum_{k=1}^n(-1)^k\dfrac{\binom{n+1}{k+1}}{\binom{n+m k}{n}}\stf{n+m k}{m k}\,. 
$$ 
\label{Proposition5}
\end{Prop}  

For example, when $m=3$, we have 
$$
c_1^{(3)}=\frac{3}{2},~c_2^{(3)}=1,~c_3^{(3)}=0,~c_4^{(3)}=\frac{1}{10},~c_5^{(3)}=-\frac{1}{4},~c_6^{(3)}=\frac{16}{21},~c_7^{(3)}=-\frac{11}{4},~c_8^{(3)}=\frac{329}{30}\,. 
$$ 
 
\bigskip

Define the Cauchy-Carlitz numbers $CC_n^{(m)}$ of order $m$ by 
\begin{equation} 
\left(\frac{x}{\log_C(x)}\right)^m=\sum_{n=0}^\infty\frac{CC_n^{(m)}}{\Pi(n)}x^n\,. 
\label{carlitzcauchyhighorder} 
\end{equation}
%By applying Lemma \ref{quotientrules} (\ref{quotientrule1}) first and then Lemma \ref{productrule2} to (\ref{carlitzcauchyhighorder}), and using the identity (\ref{eq20}), 
In the rest of this section, we show that Cauchy-Carlitz numbers of higher order are also 
expressed only in terms of \(L_i\) in the same way as Theorem \ref{cc_ht}. \par
Let again $g(z)=(\log_C(z))/z$. Applying Lemma \ref{quotientrules} (\ref{quotientrule1}) with 
$f(z)=(g(z))^m$, we get, by (\ref{carlitzcauchyhighorder}), 
$$
\frac{CC_n^{(m)}}{\Pi(n)}
=
\sum_{k=1}^n (-1)^k \sum_{i_1,\ldots,i_k\geq 1 \atop i_1+\cdots+i_k=n}
\left.H^{(i_1)}\left(g^m\right)\right|_{z=0} \cdots \left.H^{(i_k)}\left(g^m\right)\right|_{z=0}\,.
$$
By applying Lemma \ref{productrule2} to the right-hand side of the equality above, 
and using the identity (\ref{coefficientg}), 
we get the following:   
\begin{Prop}  
For $n\ge 1$, 
$$
CC_n^{(m)}=\Pi(n)\sum_{k=1}^n(-1)^k\sum_{i_1,\dots,i_k\ge 1\atop i_1+\cdots+i_k=n}M^{(m)}(i_1)\cdots M^{(m)}(i_k)\,, 
$$ 
where 
$$
M^{(m)}(i)=\sum_{j_1,\ldots, j_m\geq 0%r^{j_1},\dots,r^{j_m}\ge 1
\atop r^{j_1}+\cdots+r^{j_m}=i+m}\frac{(-1)^{j_1+\cdots+j_m}}{L_{j_1}\cdots L_{j_m}}\,. 
$$ 
\label{highorderccarlitz}
\end{Prop}

\section{Stirling-Carlitz numbers}   

One of the most useful identities of Stirling numbers is the pair of inversion properties: 
\begin{align*} 
\sum_{m=k}^n(-1)^{n-m}\stf{n}{m}\sts{m}{k}&=\delta_{n,k}\,,\\
\sum_{m=k}^n(-1)^{m-k}\sts{n}{m}\stf{m}{k}&=\delta_{n,k}\,. 
\end{align*}  
Stirling-Carlitz numbers also satisfy the similar orthogonal identities.  

%Let \(n\) be a nonnegative integer with \(r\)-ary expansion 
%\[n=\sum_{i=0}^t s_i r^i,\]
%where \(0\leq s_i\leq r-1\) for any \(i\). 

\begin{theorem}
Let $n$, $k$ be nonnegative integers with $n\geq k$. Then 
\begin{align}
&\sum_{m=k}^n\stf{n}{m}_C \sts{m}{k}_C=\delta_{n,k}, 
\label{eqnarray1}
\\
&\sum_{m=k}^n\sts{n}{m}_C \stf{m}{k}_C=\delta_{n,k}.
\label{eqnarray1b}
\end{align}
\label{orthogonality}
\end{theorem}
\begin{proof}
We may assume that $k\geq 1$ because 
if $k=0$, then (\ref{eqnarray1}) and (\ref{eqnarray1b}) are easily checked 
by (\ref{eqn:aabb}) and (\ref{eqn:aacc}). 
We see that 
\begin{align*}
z^k&=\Bigl(e_C\bigl(\log_C(z)\bigr)\Bigr)^k\\
&=
\sum_{m=k}^{\infty}\sts{m}{k}_C\frac{\Pi(k)}{\Pi(m)} \bigl(\log_C(z)\bigr)^m
\\
&=
\sum_{m=k}^{\infty}\sts{m}{k}_C\frac{\Pi(k)}{\Pi(m)} 
\sum_{n=m}^{\infty}\stf{n}{m}_C\frac{\Pi(m)}{\Pi(n)}z^n\\
&=
\sum_{n=k}^{\infty}\frac{\Pi(k)}{\Pi(n)}z^n \sum_{m=k}^n\stf{n}{m}_C \sts{m}{k}_C\,, 
\end{align*}
which implies (\ref{eqnarray1}). In the same way, we deduce (\ref{eqnarray1b}).
\end{proof}

In the rest of this section we introduce more properties on Stirling-Carlitz numbers of the first and second kind $\stf{n}{m}_C$ and $\sts{n}{m}_C$ in the case where 
$n$ and $m$ satisfy certain conditions on the sum of digits. 
When a nonnegative integer $i$ is expressed as $r$-ary expansion (\ref{r-expansion}), 
we write the sum of digits of $i$ by 
$$
\lambda(i):=\sum_{j=0}^m c_j\in\mathbb{Z}.
$$ 
First we show that if $\lambda(n)=\lambda(m)=1$, 
then the Stirling-Carlitz numbers are expressed as 
certain products of the terms of \(D_i\) and \(L_j\). 
%First we consider the case of $\lambda(n)=\lambda(m)=1$. 

\begin{Prop}
Let $a$ and $b$ be nonnegative integers with $a\geq b$. 
Then we have the following: 
\begin{align}
& \stf{r^a}{r^b}_C
=
\frac{D_a}{D_b}\cdot\frac{(-1)^{a-b}}{L_{a-b}^{r^b}}\,, 
\label{stab}
\\
& \sts{r^a}{r^b}_C
=\frac{D_a}{D_b}\cdot\frac{1}{D_{a-b}^{r^b}}\,.
\label{stac}
\end{align}
\end{Prop}
\begin{proof}
Note for any $i\geq 0$ that 
$(-1)^{r^i}=-1$. 
In fact, if $r$ is even, then $-1=1$ because the characteristic of $\mathbb{F}_r$ is $2$. 
Since 
\begin{align*}
\bigl(\log_C(z)\bigr)^{r^b}
&=
\sum_{i=0}^{\infty}\stf{r^i}{r^b}_C \frac{D_b}{D_i} z^{r^i}\\
&=
\sum_{i=0}^{\infty}\frac{(-1)^i}{L_i^{r^b}} z^{r^{i+b}}
=
\sum_{i=b}^{\infty}\frac{(-1)^{i-b}}{L_{i-b}^{r^b}} z^{r^{i}},
\end{align*}
we obtain (\ref{stab}). Similarly, calculating the coefficients of $\bigl(e_C(z)\bigr)^{r^b}$, we see 
(\ref{stac}). 
\end{proof} 

Next, we show that if $\lambda(n)>\lambda(m)$, then the Stirling-Carlitz numbers vanish. 
%Next, we consider the case of $\lambda(n)>\lambda(m)$. 

\begin{Prop}
Let $n$ and $m$ be positive integers with $\lambda(n)>\lambda(m)$. Then 
$$
\stf{n}{m}_C=\sts{n}{m}_C=0\,.
$$
\label{propst}
\end{Prop}
\begin{proof}
Let $s:=\lambda(m)$. Then $(\log_C (z))^m$ is written as 
$$
\bigl(\log_C(z)\bigr)^m=
\sum_{j=1}^{\infty}\stf{j}{m}_C\frac{\Pi(m)}{\Pi(j)} z^j
=\prod_{i=1}^s \bigl(\log_C(z)\bigr)^{r^{a(i)}}\,,
$$
where $a(1),\ldots,a(s)$ are nonnegative integers. 
For example, $a(1)=\cdots=a(c_0)=0$, 
$a(c_0+1)=\cdots=a(c_1)=1$, $\dots$, 
$a(c_{m-1}+1)=\cdots=a(c_m)=m$. 
For any $j\geq 1$, the coefficient of $z^j$ is not zero only if $\lambda(j)\leq \lambda(m)$. 
Thus, we get 
$$
\stf{n}{m}_C=0.
$$
Similarly, we see 
$$
\sts{n}{m}_C=0.
$$
\end{proof}

For example, consider the case of $r=3$, $n=8$ and $m=2$. Using $\lambda(8)=4, \lambda(2)=2$, 
we see $$\stf{8}{2}_C=\sts{8}{2}_C=0$$ by Proposition \ref{propst}.

\section{Some properties of Cauchy-Carlitz numbers}  

It is known that poly-Cauchy numbers $\mathfrak c_n^{(k)}$ satisfy 
$$
\sum_{m=0}^n\sts{n}{m}\mathfrak c_m^{(k)}=\frac{1}{(n+1)^k}
$$ 
(\cite[Theorem 3]{Ko1}).  If $k=1$, this identity is the same as that in \cite[Theorem 2.3]{MSV}.  For Cauchy-Carlitz numbers, we obtain an analogous identity.

\begin{theorem}  
For a nonnegative integer $n$, we have 
$$
\sum_{m=0}^{n}\sts{n}{m}_C CC_m
=\begin{cases}\frac{1}{L_j}&\text{if $n=r^j-1$}\,,\\
0&\text{otherwise}\,.
\end{cases}   
$$ 
\label{th33}
\end{theorem}  
\begin{proof}  
By Theorem \ref{th10} and Theorem \ref{orthogonality}, we have 
\begin{align*}  
\sum_{m=0}^{n}\sts{n}{m}_C CC_m&=
\sum_{m=0}^n\sts{n}{m}_C\sum_{j=0}^\infty\frac{1}{L_j}\stf{m}{r^j-1}_C\\
&=\sum_{j=0}^\infty\frac{1}{L_j}\sum_{m=0}^n\sts{n}{m}_C\stf{m}{r^j-1}_C\\
&=\sum_{j=0}^\infty\frac{1}{L_j}\delta_{n,r^j-1}\\
&=\begin{cases}\frac{1}{L_j}&\text{if $n=r^j-1$}\,,\\
0&\text{otherwise}\,.
\end{cases}   
\end{align*} 
\end{proof}

It is known that 
$$
\frac{1}{n!}\sum_{m=0}^n(-1)^m\stf{n+1}{m+1}B_m=\frac{1}{n+1}\,. 
$$ 
Similarly, we have the following.  

\begin{theorem}  
For a nonnegative integer $n$, we have 
$$
\sum_{m=0}^{n}\stf{n}{m}_C BC_m
=\begin{cases}\frac{(-1)^j D_j}{L_j^2}&\text{if $n=r^j-1$}\,,\\
0&\text{otherwise}\,.
\end{cases}   
$$ 
\label{th44}
\end{theorem}

There are alternating expressions between poly-Bernoulli numbers and poly-Cauchy numbers (\cite{Ko1,KL}). When $k=1$, they are reduced to the relations between classical Bernoulli numbers and classical Cauchy numbers.  
\begin{align*} 
B_n^{(k)}&=\sum_{l=0}^n\sum_{m=0}^n(-1)^{n-m}m!\sts{n}{m}\sts{m}{l}\mathfrak c_l^{(k)}\,,\\
\mathfrak c_n^{(k)}&=\sum_{l=0}^n\sum_{m=0}^n\frac{(-1)^{n-m}}{m!}\stf{n}{m}\stf{m}{l}B_l^{(k)}\,. 
\end{align*}  

As analogues, we have the following.  

\begin{theorem} 
\begin{align*} 
BC_n&=\sum_{l\ge 0}\sum_{m=r^j-1\ge 0}(-1)^{j}\Pi(m)\sts{n}{m}_C\sts{m}{l}_C CC_l\,,\\
CC_n&=\sum_{l\ge 0}\sum_{m=r^j-1\ge 0}\frac{(-1)^{j}}{\Pi(m)}\stf{n}{m}_C\stf{m}{l}_C BC_l\,. 
\end{align*} 
\label{th55} 
\end{theorem} 

\begin{proof} 
By Theorem \ref{th33} and Theorem \ref{th20}, we have 
\begin{align*}
&\sum_{l\ge 0}\sum_{m=r^j-1\ge 0}(-1)^{j}\Pi(m)\sts{n}{m}_C\sts{m}{l}_C CC_l\\
&=\sum_{j=0}^\infty(-1)^j\Pi(r^j-1)\sts{n}{r^j-1}_C\frac{1}{L_j}\\
&=\sum_{j=0}^\infty\frac{(-1)^j D_j}{L_j^2}\sts{n}{r^j-1}_C\\
&=BC_n\,. 
\end{align*} 
By Theorem \ref{th44} and Theorem \ref{th10}, we have 
\begin{align*} 
&\sum_{l\ge 0}\sum_{m=r^j-1\ge 0}\frac{(-1)^{j}}{\Pi(m)}\stf{n}{m}_C\stf{m}{l}_C BC_l\\
&=\sum_{j=0}^\infty\Pi(r^j-1)\stf{n}{r^j-1}_C\frac{(-1)^j D_j}{L_j^2}
=\sum_{j=0}^\infty\frac{1}{L_j}\stf{n}{r^j-1}_C\\
&=CC_n\,. 
\end{align*} 
\end{proof} 

\section*{Acknowledgements} 

The first author is supported by JSPS KAKENHI Grant Number 15K17505.  
The second author is in part supported by the grant of Wuhan University and by the grant of Hubei Provincial Experts Program.

\end{document}